\documentclass[10pt,leqno]{article}
\usepackage{amssymb,amsbsy,amsmath,amsfonts,amssymb,amscd, mathrsfs}

\usepackage[english]{babel}
\usepackage[T1]{fontenc}
\usepackage{indentfirst}

\usepackage{ifpdf}
\ifpdf
	\usepackage[pdftex]{graphicx}
\else
	\usepackage[dvips]{graphicx}
\fi

\makeatletter
\@addtoreset{equation}{section}
\makeatother

\newtheorem{statement}{}[section]
\newtheorem{theoreme}[statement]{Theorem}
\newtheorem{lemme}[statement]{Lemma}

\newtheorem{proposition}[statement]{Proposition}

\newtheorem{corollaire}[statement]{Corollary}

\newcommand\C{\mathbb C}

\newcommand\T{\mathbb T}
\newcommand\D{\mathbb D}

\newcommand\e{{\rm e}}
\renewcommand\P{\mathbb P}
\newcommand\eps{\varepsilon}
\newcommand\ind{{\rm 1\kern-.30em I}}
\newcommand\qed{\hfill $\square$}
\renewcommand \Re{{\mathfrak R}{\rm e}\,}
\renewcommand \Im{{\mathfrak I}{\rm m}\,}
\let\phi=\varphi

\title{\bf Some revisited results about composition operators on Hardy spaces}
\author{\it Pascal Lef\`evre, Daniel Li,\\ \it Herv\'e Queff\'elec, Luis Rodr{\'\i}guez-Piazza}

\date{\footnotesize \today}

\begin{document}

\maketitle

\noindent{\bf Abstract.} \emph{We generalize, on one hand, some results known for composition operators on 
Hardy spaces to the case of Hardy-Orlicz spaces $H^\Psi$: construction of a ``slow'' Blaschke product giving a 
non-compact composition operator on $H^\Psi$; construction of a surjective symbol whose composition operator 
is compact on $H^\Psi$ and, moreover, is in all the Schatten classes $S_p (H^2)$, $p > 0$. On the other hand, we 
revisit the classical case of composition operators on $H^2$, giving first a new, and simplier, characterization of 
closed range composition operators, and then showing directly the equivalence of the two characterizations of 
membership in the Schatten classes of Luecking and Luecking and Zhu.} 
\medskip

\medskip



\section{Introduction}

The study of composition operators on Hardy spaces is now a classical subject (see \cite{Shap-livre}, 
\cite{Co-McC} for example). In \cite{CompOrli} (see also \cite{LLQR-CRAS}), we considered a more general setting 
and studied composition operators on Hardy-Orlicz spaces; we gave there a characterization of their compactness in 
terms of the Carleson function of their symbol (and in terms of the Nevanlinna counting function in \cite{LLQR-N}). 
This work was continued in \cite{JMAA}: we compared the compactness on Hardy spaces versus the  compactness 
on Hardy-Orlicz spaces. For instance, we showed that there is, for every $1 \leq p < \infty$, an Orlicz function 
$\Psi$ such that $H^{p + \eps} \subseteq H^\Psi \subseteq H^p$ for every $\eps > 0$, and a composition 
operator $C_\phi$ such that $C_\phi$ is compact on $H^p$ and $H^{p + \eps}$, but which is not compact on 
$H^\Psi$.\par
We carry on this study in the present work. In a first part (Section~\ref{section slow Blaschke} and 
Section~\ref{surjective symbol}), we shall improve, and extend to the Hardy-Orlicz case, results known 
for Hardy spaces; in a second part (Section~\ref{section closed range} and Section~\ref{Schatten}), we shall give 
new lights on some results concerning Hardy spaces. More precisely, the content of this paper is as following.\par
\smallskip

B. McCluer and J. Shapiro (\cite{McCluer-Shapiro}, Theorem~3.10; see also \cite{Shap-livre}, \S~3.2) proved that, 
when their symbol $\phi$ is finitely-valent, compactness of composition operators $C_\phi$ on the Hardy space 
$H^2$ can be characterized by the behaviour of the modulus of $\phi$ near the frontier of $\D$: compactness is 
equivalent to $1 - |z| = 0\, \big( 1 - |\phi (z)| \big)$ as $|z| \to 1$, but that is not equivalent in general 
(\cite{McCluer-Shapiro}, Example~3.8; see also \cite{Shap-livre}, \S~10.2). In \cite{LLQR-N}, Theorem~5.3, we 
gave such a characterization for composition operators, with finitely-valent symbol, on Hardy-Orlicz spaces.  
In Section~\ref{section slow Blaschke}, we construct a ``slow'' Blaschke product (generalizing 
\cite{Shap-livre}, \S~10.2 and \cite{CompOrli}, Proposition~5.5) showing that this condition is not sufficient in 
general.\par
\smallskip

In Section~\ref{surjective symbol}, we construct a compact composition operator 
$C_\phi \colon H^\Psi \to H^\Psi$ with surjective symbol $\phi$ and such that $C_\phi \colon H^2 \to H^2$ 
is in all the Schatten classes $S_p (H^2)$, $p > 0$. This generalizes and improves a result of 
B. McCluer and J. Shapiro (\cite{McCluer-Shapiro}, Example~3.12; see also the survey \cite{Queff}, \S~2). \par
\smallskip

In Section~\ref{section closed range}, we give a characterization of composition operators 
$C_\phi \colon H^p \to H^p$, $1 \leq p < \infty$, with a closed range, simpler than 
the former ones (see \cite{Cima} and \cite{Nina}).\par
\smallskip

Finally, based on the main result of \cite{LLQR-N}, we show directly, in Section~\ref{Schatten}, the equivalence of 
Luecking's and Luecking-Zhu's criteria  (\cite{Luecking}, \cite{Luecking-Zhu}) for the membership of 
$C_\phi \colon H^2 \to H^2$ in the Schatten classes.\par
\bigskip

\noindent{\bf Acknowledgement.} Part of this work was made during the fourth-named 
author visited the University of Lille 1 and the University of Artois (Lens) in June 2009. This 
fourth-named author is partially supported by a Spanish research project MTM2006-05622.


\section{Notation}

The open unit disk is denoted by $\D = \{z \in \C\,;\ |z| < 1\}$ and its boundary, the unit circle, by 
$\T = \{z \in \C\,;\ |z| = 1\}$. The normalized Lebesgue measure $dt/2\pi$ on $\T$ is denoted by $m$. The 
normalized area measure $dx\,dy/\pi$ is denoted by $A$.\par 
The Hardy space $H^1$ is the space of analytic functions $f \colon \D \to \C$ such that 
$\sup_{r < 1} \int_0^{2\pi} | f (r \e^{i\theta})|\, d\theta < \infty$. Every $f \in H^1$ has almost everywhere 
boundary values on $\T$, which are denoted by $f^\ast$.\par
An Orlicz function is a convex nondecreasing function $\Psi \colon [0, \infty) \to [0, \infty)$ such that 
$\Psi (0) =0$ and $\Psi (\infty) = \infty$. If $\mu$ is a positive measure on some measurable space $S$, the 
Orlicz space $L^\Psi (\mu)$ is the set of all (classes of) measurable functions $f \colon S \to \C$ such that 
$\int_S \Psi (|f|/ C)\, d\mu < \infty$ for some $C > 0$; the norm $\| f \|_\Psi$ is defined as the infimum of the 
positive numbers $C$ for which $\int_S \Psi (|f|/ C)\, d\mu \leq 1$.\par
The Hardy-Orlicz space $H^\Psi$ is the linear subspace of $f \in H^1$ such that 
$f^\ast \in L^\Psi (m)$ (see \cite{CompOrli}).\par\smallskip

Every analytic self-map $\phi \colon \D \to \D$ defines a bounded composition operator 
$C_\phi \colon f\in H^\Psi \mapsto f \circ \phi \in H^\Psi$ (see \cite{CompOrli}).\par
\smallskip

For every $\xi \in \T$ and $0 < h < 1$, the Carleson window is the set 
$W (\xi, h) = \{z\in \D\,;\ |z| \geq 1 - h \ \text{and} \ |\arg (z\,\bar{\xi} | \leq h\}$. The Carleson function 
$\rho_\phi$ of the analytic self-map $\phi \colon \D \to \D$ is defined, for $0  < h < 1$, by:
\begin{displaymath}
\rho_\phi (h) 
= \sup_{\xi \in \T} m \big( \{\e^{i \theta}\in \T \,;\ \phi^\ast (\e^{i \theta}) \in W (\xi, h) \} \big)\,.
\end{displaymath}
Alternatively, $\rho_\phi (h) = \sup_{\xi \in \T} m_\phi [W (\xi, h)]$, where $m_\phi$ is the pull-back measure 
of $m$ by $\phi$. We shall also use, instead of $W (\xi, h)$, the set 
$S (\xi, h) =\{ z\in \D\,;\ |z - \xi | \leq h\}$, which has an equivalent size.\par
\smallskip

The Nevanlinna counting function $N_\phi$ is defined, for $w \in \phi(\D) \setminus \{\phi (0) \}$, by
\begin{displaymath}
N_\phi (w) = \sum_{\phi (z) = w} \log \frac{1}{|z|}\,,
\end{displaymath}
each term $\log \frac{1}{|z|}$ being repeated according to the multiplicity of $z$, and $N_\phi (w) = 0$ 
for the other $w \in \D$.


\section{Slow Blaschke products}\label{section slow Blaschke}

B. McCluer and J. Shapiro (\cite{McCluer-Shapiro}, Theorem~3.10; see also \cite{Shap-livre}, \S~3.2) proved that, 
when $\phi$ is finitely-valent (meaning that, for some $s \geq 1$, the equation $\phi (z) = w$ has at most $s$ 
solutions), the composition operators $C_\phi \colon H^p \to H^p$ is compact, $1 \leq p < \infty$, if and only 
if $\phi$ has an angular derivative at no point of $\T$; that means that:
\begin{equation}\label{angular derivative}
\lim_{|z| \to 1} \frac {1 - |z|} {1 - |\phi (z)|} = 0\,.
\end{equation}

In \cite{LLQR-N}, Theorem~5.3, we generalized this result to Hardy-Orlicz spaces and proved that if $\phi$ is 
finitely-valent, the composition operator $C_\phi \colon H^\Psi \to \colon H^\Psi$ is compact if and only if:
\begin{equation}\label{condition compacite Nevanlinna}
\lim_{|z| \to 1} \frac{\Psi^{- 1} \bigg[\displaystyle  \frac{1}{1 - |\phi (z) |} \bigg]} 
{\Psi^{-1} \bigg[\displaystyle \frac{1}{1 - |z|} \bigg]} = 0\,.
\end{equation}
Without the assumption that $\phi$ is finitely-valent, condition~\eqref{condition compacite Nevanlinna} is 
no longer sufficient to ensure the compactness of $C_\phi \colon H^\Psi \to H^\Psi$. Indeed, we are going to 
construct a Blaschke product satisfying \eqref{condition compacite Nevanlinna}, but whose associated 
composition operator is of course not compact on $H^\Psi$, as this is the case for every inner function. A  
Blaschke product satisfying \eqref{angular derivative} is constructed in \cite{Shap-livre}, \S~10.2; 
that construction uses Frostman's Theorem. Our construction, which is more general, is entirely elementary.

\begin{theoreme}\label{slow Blaschke}
Let $\delta\colon (0,1) \to (0, 1/2]$ be any function such that $\lim\limits_{t \to 0} \delta (t) = 0$. 
Then, there exists a Blaschke product $B$ such that:
\begin{equation}\label{condition slow Blaschke}
\qquad \qquad 1 - |B(z)| \geq \delta(1 - |z|), \qquad \text{for all $z\in \D$.}
\end{equation}
\end{theoreme}

\begin{corollaire}
For every Orlicz function $\Psi$ there exists a Blaschke product $B$ which satisfies:
\begin{displaymath}
\lim_{|z| \to 1} \frac{\Psi^{-1} \bigg[\displaystyle  \frac{1}{1 - |B (z) |} \bigg]} 
{\Psi^{-1} \bigg[\displaystyle \frac{1}{1 - |z|} \bigg]} = 0\,.
\end{displaymath}
though the composition operator $C_B \colon H^\Psi \to H^\Psi$ is not compact.
\end{corollaire}

\noindent{\bf Proof.} $C_B$ is not compact since every compact composition operator should satisfy 
$| \phi^\ast | < 1$ \emph{a.e.} (see \cite{CompOrli}, Lemma~4.8). It suffices then to chose 
$ \delta (t) = 1 / \Psi \big( \sqrt{ \Psi^{-1} (1/t) }\big)$, which satisfies the hypothesis of 
Theorem~\ref{slow Blaschke}. Moreover:
\begin{displaymath}
\frac{\Psi^{- 1} \big(1/ \delta (t) \big)}{\Psi^{- 1} (1/t)} = \frac{1}{\sqrt{\Psi^{- 1} (1/t)}} \hskip 2pt 
\mathop{\longrightarrow}_{t \to 0} 0\,,
\end{displaymath}
and condition \eqref{condition slow Blaschke} gives the result. \qed
\medskip

\noindent{\bf Proof of Theorem~\ref{slow Blaschke}.} 
We shall essentially construct our Blaschke product $B$ as an infinite product of finite Blaschke products
\begin{equation*}
\prod_n B_n \,,
\end{equation*}
where each finite Blaschke product $B_n$ has $p_n$ zeros equidistributed in the circumference of radius $r_n$. 
That is, we will have, writing $\theta_k = 2\pi k/p_n$ and $z_k = r_n\, \e^{i \theta_k}$, for 
$k = 1, 2, \dots, p_n$:
\begin{equation}\label{finito}
B_n (z) =\prod_{k = 1}^{p_n} 
\frac{|z_k|}{z_k} \frac{z_k - z}{1 - \overline{z}_k z}
=\prod_{k = 1}^{p_n} \frac{r_n - \e^{-i \theta_k} z}{1 - r_n \e^{-i \theta_k} z} \,\cdot
\end{equation}

We shall need the following estimate for the finite Blaschke product in \eqref{finito}.

\begin{lemme}\label{lemme Blaschke}
Let $p\in \mathbb{N}$, and $0 < r < 1$. Consider the finite Blaschke product
\begin{equation}
G (z) = \prod_{k=1}^p \frac{r - \e^{-i \theta_k} z}{1 - r \e^{-i\theta_k} z} \,\raise 1pt \hbox{,}
\end{equation}
where $\theta_k =\frac{2k\pi}{p}$, for $k = 1, 2,\ldots, p$.\par

$(a)$ Then, for every $z\in \mathbb{D}$ with $|z|=r$,
\begin{equation}\label{hyperbol}
|G (z)| \le \frac{2r^p}{1 + r^{2p}} = 1 - \frac{(1 - r^p)^2}{1 + r^{2p}} \,\cdot
\end{equation}
\par

$(b)$ If besides we have $p\,h \le 1/2$, where $h = 1 - r$, we also have, for every $z\in \mathbb{D}$ 
with $|z| = r$,
\begin{equation}
|G (z)| \le 1 - \frac{(p\,h)^2}{2 \e} \,\cdot
\end{equation}
\end{lemme}

Let us continue the proof of the theorem. Define $\chi \colon (0, 1) \to (0, 1]$ by:
\begin{equation}\label{definition chi}
\chi (x) = \sup_{t \leq x} \big[\max\{ 2\delta (t), \sqrt t\} \big] \,.
\end{equation}
Then $\chi$ is non-decreasing, $\lim_{x \to 0} \chi (x) = 0$ and $\lim_{x \to 1} \chi (x) = 1$. 
We can find a decreasing sequence $(h_n)_{n \geq 0}$ of point $h_n\in (0, 1)$, such that 
$\chi (h_n) \leq 2^{-n}$. This sequence converges to $0$; in fact, $\sqrt{h_n} \leq \chi (h_n) \leq 2^{-n}$, 
by \eqref{definition chi},  and hence:
\begin{equation}
h_n \leq 2^{- 2n} \,.
\end{equation}

We now define, for every $n\in \mathbb{N}$, a positive integer $p_n$, by:
\begin{equation}\label{defino}
p_n = \min \{ p\in \mathbb{N}\,; \  \frac{p^2h_n^2}{2\e} > 2^{-n}\  \}.
\end{equation}

We have $p_n > 1$ because $h_n^2/2 \e < h_n^2 \le 2^{-4n}$. So, for every $n$, we have 
$4(p_n - 1)^2 \ge p_n^2$, and then:
\begin{equation}\label{tres}
4\cdot 2^{-n} \ge \frac {4(p_n-1)^2h_n^2}{2 \e} \ge \frac{p_n^2h_n^2}{2 \e} \,\cdot
\end{equation}
This yields, for $n\ge 7$, that $(p_n h_n)^2 \le 8 \e \, 2^{-n} \le 1/4$. 
Therefore $p_n h_n \le 1/2$, and we can use the estimate in part $(b)$ of Lemma~\ref{lemme Blaschke}.
\smallskip

Now, for $n \ge 7$, let $B_n$ be the finite Blaschke product defined by \eqref{finito}, where $r_n = 1 - h_n$. 
Using $(b)$ in Lemma~\ref{lemme Blaschke}, the Maximum Modulus Principle and the definition of $p_n$ 
in \eqref{defino}, we have:
\begin{equation}\label{cuatro}
\qquad \qquad |B_n (z)| \le 1 - \frac{p_n^2 h_n^2}{2 \e} < 1 - 2^{-n}, 
\quad \hbox{for $|z| \le r_n$.}
\end{equation}

Consider then the Blaschke product $D$ defined by:
\begin{equation}
D (z) =\prod_{n=7}^\infty B_n (z).
\end{equation}
This product is convergent since, by \eqref{tres}, we have:
\begin{equation*}
\sum p_n (1 - r_n) = \sum p_n h_n \le \sum \sqrt{8 \e \,2^{-n}} < +\infty\,.
\end{equation*}

Finally, take $N\in \mathbb{N}$ big enough to have $r_6^N < 1/2$, 
and define:
\begin{equation}
B (z)= z^N \,D(z).
\end{equation}

Thus $B$ is a Blaschke product, and, if $|z| \le r_6$, we have, since $\delta (t) \leq 1/2$:
\begin{equation}\label{menorR}
|B(z)| \le | z^N| \le r_6^N < 1/2 \le 1 - \delta(1 - |z|).
\end{equation}

If $1 > |z| > r_6$, there exists $k \ge 7$ such that $r_k \ge |z| > r_{k - 1}$. Therefore, thanks to \eqref{cuatro},
\begin{equation}\label{mayorRa}
|B (z)| \le |D (z)| \le |B_k (z)| \le 1 - 2^{-k}.
\end{equation}
On the other hand $r_k \ge |z| > r_{k - 1}$ implies $h_k \le 1 - |z| < h_{k - 1}$, and so:
\begin{equation}\label{mayorRb}
\delta (1 - |z|) \le \frac{1}{2} \chi(1 - |z|) \le  \frac{1}{2} \chi(h_{k-1}) \leq 2^{- k} \,.
\end{equation}

Combining \eqref{mayorRa} and \eqref{mayorRb} we get $|B (z)| \le 1 -\delta (1 - |z|)$, when 
$1 > |z| > r_6$. From this and \eqref{menorR}, Theorem~\ref{slow Blaschke} follows.
\qed
\medskip

\noindent{\bf Proof of Lemma~\ref{lemme Blaschke}.} It is obvious that, for all  $a, z \in \mathbb{C}$, 
\begin{displaymath}
\qquad \qquad \prod_{k = 1}^p (z - a \e^{i\theta_k}) = z^p - a^p \, .
\end{displaymath}

Using this we have:
\begin{equation}
G (z) =\prod_{k=1}^p \frac{r - \e^{-i\theta_k} z}{1 - r \e^{-i\theta_k}z} =
\prod_{k=1}^p \frac{z - r \e^{i\theta_k}}{r z - \e^{i\theta_k}} =
\frac{z^p - r^p}{(rz)^p - 1} \,\cdot
\end{equation}
Now, if $|z| = r$, we can write $z^p = r^p u$, for some $u$ with $|u| = 1$. Then $|G (z)| = |T (u)|$, where $T$ 
is the Moebius transformation
\begin{displaymath}
T (u) = \frac{r^p(u - 1)}{r^{2p} u - 1} \,\cdot
\end{displaymath}
This transformation $T$ maps the unit circle $\partial\mathbb{D}$ onto a circumference $C$. 
As $T$ maps the extended real line $\mathbb{R}_\infty$ to itself, and $\partial\mathbb{D}$ is orthogonal 
to $\mathbb{R}_\infty$ at the intersection points $1$ and $-1$, $C$ is the circumference orthogonal to 
$\mathbb{R}_\infty$ crossing through the points $T (1) = 0$ and $T (- 1) =\alpha$. It is easy to see that 
$|w|\le |\alpha|$, for every $w\in C$; consequently:
\begin{displaymath}
|G (z)| \le \sup_{u \in \partial\mathbb{D} } |T(u)| = |T (-1)| = \frac{2r^p}{1 + r^{2p}} \,\cdot
\end{displaymath}
This finishes the proof of the statement $(a)$.\par
\smallskip

To prove part (b), observe that, $1 + r^{2p} \le 2$, and so, for $|z|=r$,
\begin{equation}\label{ultima}
|G (z)| \le 1 - \frac{(1 - r^p)^2}{1 + r^{2p}} \le 1 - \frac{(1 - r^p)^2}{2} \,\cdot
\end{equation}
Remember that $r = 1 - h$, so $r \le \e^{-h}$, and $r^p \le \e^{-ph}$. 
Thus $1 - r^p\ge 1- \e^{-ph}$. Now, if $x \in [0, 1/2]$, we have, by the Mean Value theorem:
\begin{displaymath}
1 - \e^{- x}  \ge \frac{x}{\sqrt \e} \,\cdot
\end{displaymath}
Since $p\,h \le 1/2$, we can apply this last estimate to \eqref{ultima} to get, as promised,
\begin{displaymath}
|G (z)| \le 1 -\frac{(1 - \e^{-ph})^2}{2} \le 1 - \frac{p^2h^2}{2 \e} \, \raise 1pt \hbox{,}
\end{displaymath}
and ending the proof of Lemma~\ref{lemme Blaschke}. \qed
\medskip

\noindent{\bf Remark.} The key point in the proof of Theorem~\ref{slow Blaschke} is the inequality 
\eqref{hyperbol} in Lemma~\ref{lemme Blaschke}. This inequality may be viewed as a consequence of the 
strong triangle inequality (applied to $a = z^p$, $b = r^p$ and $c = 0$):
\begin{equation}\label{triangle}
d (a, b) \leq \frac{d(a, c) + d (c, b)}{1 + d (a, c)\, d(c, b)}
\end{equation}
for the pseudo-hyperbolic distance $d (u, v) = \frac{|u - v|}{|1 - \bar{u} v|}$ on $\D$. Let us recall a proof 
for the convenience of the reader: by conformal invariance, we may assume that $c = 0$; then:
\begin{displaymath}
1 - [d (a, b)]^2 = \frac{(1 - |a|^2) (1 - |b|^2)}{|1 - \bar{a} b|^2}  
\geq \frac{(1 - |a|^2) (1 - |b|^2)}{(1 + |a|\,|b|)^2} = 1 - [d (|a|, - |b|)]^2 \,, 
\end{displaymath}
so that:
\begin{displaymath}
d (a, b) \leq d (|a|, -|b|) = \frac{|a| + |b|}{1 + |a|\, |b|} \,\raise 1pt \hbox{,}
\end{displaymath}
proving \eqref{triangle}, since $d (a, 0) = |a|$ and $d (0, b) = |b|$.
\par\goodbreak


\section{A compact composition operator with a surjective symbol}\label{surjective symbol}
 
A well-known result of J. H. Schwartz (\cite{Schw}, Theorem~2.8) asserts that the composition operator 
$C_\phi \colon H^\infty \to H^\infty$ is compact if and only if $\| \phi \|_\infty < 1$. In particular, the 
compactness of $C_\phi \colon H^\infty \to H^\infty$ prevents the surjectivity of $\phi$. It may be therefore 
to be expected that, the bigger $\Psi$, the more difficult it will be to obtain both the compactness of 
$C_\phi \colon H^\Psi \to H^\Psi$ and the surjectivity of $\phi$. Nevertheless, this is possible, as 
says the following theorem, and the case $H^\infty$ appears really as a singular case (corresponding to an 
``Orlicz function'' which is discontinuous and can take the value infinity).

\begin{theoreme} \label{surjective}
For every Orlicz function $\Psi$, there exists a symbol $\varphi \colon \D \to \D$ which is 
$4$-valent and surjective  and such that $C_\varphi \colon H^\Psi \to H^\Psi$ is compact. Moreover, 
$\varphi$ can be taken so as $C_\varphi \colon H^2 \to H^2$ is in all the Schatten classes $S_p (H^2)$, $p > 0$.
\end{theoreme}

In the case of $H^2$ ($\Psi (x) = x^2$), B. McCluer and J. Shapiro (\cite{McCluer-Shapiro}, Example~3.12) gave 
an example based on the Riemann mapping theorem and on the fact that, for a finitely valent symbol $\varphi$, 
we have the equivalence:
\begin{equation}\label{equiv 1}
C_\varphi \colon H^2 \to H^2 \text{ compact} \quad \Longleftrightarrow \quad 
\lim_{\vert z \vert \mathop{\to}\limits^{<} 1} \frac{1 - \vert \varphi (z) \vert}{1 - \vert z\vert} =\infty.
\end{equation}

A specific example is as follows. Take 
\begin{equation}
R = \big\{ z = x+iy \in \C \,;\ x > 0 \text{ and } \frac{1}{x} < y < \frac{1}{x} + 4\pi \big\},
\end{equation}
let $g \colon \D\to R$  be a Riemann map and set $\varphi = \e^{- g}$. Then, $\varphi$ is $2$-valent, 
$\varphi (\D) = \D^\ast$ (where $\D^\ast = \D \setminus \{0\}$), and the validity of \eqref{equiv 1} is tested 
through the use of the  Julia-Carath\'eodory theorem (see \cite{Queff} for details). To get a fully surjective mapping 
$\varphi_1$, just compose $\varphi$ with the square of a Blaschke product: 
\begin{displaymath}
\qquad \qquad \varphi_{1}(z) =B \circ \varphi, \qquad \text{with}  \quad
B(z) = \Big(\frac{z - \alpha}{1 - \overline{\alpha} z} \Big)^2, \quad  \alpha\in D^\ast = \D \setminus \{0\}
\end{displaymath}
(note that $B (0) = B (2\alpha/ 1 +|\alpha|^2)$. Since $C_{\varphi_1} = C_\varphi \circ C_B$, we see that 
$C_{\varphi_1}$ is compact as well and we are done.\par \smallskip

Here, we can no longer rely on the Julia-Carath\'eodory theorem. But we shall use the following necessary and 
sufficient condition, in terms of the maximal Carleson function $\rho_{\varphi}$, which is valid for any 
symbol, finitely-valent or not (see \cite{CompOrli}, Theorem~4.18 -- or \cite{LLQR-CRAS}, Th\'eor\`eme~4.2, 
where a different, but equivalent, formulation is given):
\begin{equation}\label{equiv 2}
C_\varphi \colon H^\Psi\to H^\Psi \text{ compact} \quad \Longleftrightarrow \quad 
\lim_{ h \mathop{\to}\limits^{>} 0} \frac{\Psi^{- 1} (1/ h)}{\Psi^{- 1} \big(1/ \rho_{\varphi}(h) \big)} = 0\,.
\end{equation}
\medskip

For the sequel, we shall set:
\begin{equation}
\Delta (h) = \frac{\Psi^{- 1} (1/ h)}{\Psi^{- 1} \big(1/ \rho_{\varphi}(h) \big)} \, \cdot
\end{equation}
\par\medskip

Our strategy will be to elaborate on the previous example to produce a (nearly) surjective $\varphi$ such that 
$\rho_{\varphi} (h)$ is very small (depending on $\Psi$) for small $h$. The tool will be the notion of 
harmonic measure for certain open sets of the extended plane $\hat \C = \C \cup \{\infty\}$, called 
\emph{hyperbolic} (see \cite{Conway}, Definition~19.9.3); for example, every conformal image of $\D$ is 
hyperbolic (see \cite{Conway}, Proposition~19.9.2 (d) and Theorem~19.9.7). If $G$ is a hyperbolic domain and 
$a\in G$, the \emph{harmonic measure} of $G$ at $a$ is the probability measure $\omega_{G}(a, \, .\,)$ 
supported by ${\partial G}$ (here, and throughout the rest of this section, boundaries and closures will be taken 
in $\hat \C$) such that: 
\begin{displaymath}
u (a) =\int_{{\partial G}} u(z) \,d\omega_{G}(a, z) 
\end{displaymath}
for each bounded and continuous function $u$ on $\overline G$, which is harmonic in $G$ 
(see \cite{Conway}, Definition~21.1.3). The harmonic measure at $a$ of a Borel set $A\subseteq {\partial G}$ 
will be denoted by $\omega_{G} (a, A)$. Clearly, 
\begin{displaymath}
\omega_{\D}(0, \,.\,) = m, 
\end{displaymath}
the Haar measure (\emph{i.e.} normalized Lebesgue measure) of ${\partial \D}$.\par\smallskip

R. Nevanlinna (see \cite{Conway}, Proposition~21.1.6) showed that harmonic measures share a 
\emph{conformal invariance property}. Namely, assume that $G$ is a simply connected domain, in which the 
Dirichlet problem can be solved (a \emph{Dirichlet domain}), and  $\tau \colon \overline \D \to \overline G$ 
is a continuous function which maps conformally $\D$ onto $G$; then $\tau$ maps $\partial \D$ onto 
$\partial G$, and, if $\tau (0) = a$: 
\begin{equation}\label{equiv 3}
\omega_G (a, A)  = m \big( \tau^{- 1} (A) \big)  
\end{equation}
for every Borel set $A \subseteq \partial G$. This explains why harmonic measures enter the matter 
when we consider composition operators $C_{\varphi}$: such an operator induces a map 
$H^\Psi \to L^\Psi (m_\phi)$, where $m_{\varphi} =\varphi^*(m)$ appears as an image measure of $m$, 
as it happens for the harmonic measure of $G$ at $a$ in \eqref{equiv 3}.
\par\smallskip

A useful alternative way of defining the harmonic measure, due to S. Kakutani, and completed by J. Doob (see 
\cite{Stroock}, page 454, and \cite{Garnett}, Appendix~F, page~477) is the following: Let $(B_t)_{t > 0}$ be the 
$2$-dimensional Brownian motion starting at $a \in G$ (\emph{i.e.} $B_0 = a$), and $\tau$ be the stopping time 
defined by: 
\begin{equation}
\tau =\inf \{t > 0 \,;\ B_{t} \notin G \}\,; 
\end{equation}
we have:
\begin{equation}\label{proba}
\omega_{G} (a, A) =\P_{a}(B_{\tau}\in A), 
\end{equation}
\emph{i.e.} the harmonic measure of $A$ at $a$ is the probability that the Brownian motion starting at $a$ 
exits from $G$ through the Borel set $A\subseteq {\partial G}$. The following lemma will be basic for 
the construction of our example. We shall provide two proofs, the second one being more illuminating.

\begin{lemme}[Hole principle] \label{hole principal}
Let $G_0$ and $G_1$ be two hyperbolic open sets and $H \subseteq {\partial G_0}$ a Borel set such that  
\begin{displaymath}
G_0\subseteq G_1 \quad \text{and} \quad  {\partial G_0} \subseteq {\partial G_1} \cup H .
\end{displaymath}
Then, for every $a\in G_0$, we have the following inequality: 
\begin{equation}\label{inegalite trou}
\omega_{G_1} (a, {\partial G_1} \setminus {\partial G_0}) \leq \omega_{G_0} (a, H).
\end{equation}
\end{lemme}

\noindent{\bf Proof 1.} From \cite{Conway}, Corollary~21.1.14, with $\Delta = \partial G_0 \cap \partial G_1$, 
one has $\omega_{G_0} (a, \Delta) \leq \omega_{G_1} (a, \Delta)$. But 
$\partial G_1 \setminus \Delta = \partial G_1 \setminus \partial G_0$, and hence, since harmonic measures 
are probability measures,  
\begin{displaymath}
\omega_{G_1} (a, \partial G_1 \setminus \partial G_0) = \omega_{G_1} (a, \partial G_1 \setminus \Delta) 
= 1 - \omega_{G_1} (a, \Delta) \leq 1 - \omega_{G_0} (a, \Delta); 
\end{displaymath}
we get the result since 
$\partial G_0 = H \cup \Delta$, which implies $1 \leq \omega_{G_0} (a, H) + \omega_{G_0} (\Delta)$. \qed
\par
\medskip

\noindent{\bf Proof 2.} Let us define 
\begin{equation}
\tau_0 =\inf\{t  > 0 \,;\ B_t \notin G_0\},\quad \tau_1 = \inf\{t > 0 \,;\ B_t \notin G_1\}
\end{equation}
and
\begin{equation}
E = \{B_{\tau_1} \in {\partial}G_1 \setminus {\partial}G_0\}, \qquad F =\{B_{\tau_0} \in H\}.
\end{equation}

Inequality~\eqref{inegalite trou} amounts to proving that $\P_{a} (E)\leq \P_{a} (F)$, which will follow from 
the inclusion $E \subseteq F$. Suppose that the event $E$ holds. Since $G_0 \subseteq G_1$, one has 
$\tau_0 \leq \tau_1$. The Brownian path $(B_s)_{0\leq s\leq \tau_1}$ being continuous with  
$B_0 = a \in G_0$, one has $B_{\tau_0} \in {\partial G_0} \subseteq {\partial G_1} \cup H$. If we had 
$B_{\tau_0} \in {\partial G_1}$, we should have $B_{\tau_0} \notin G_1$, since $G_1$ is open, and 
hence $\tau_0 =\tau_1$, since we know that $\tau_0 \leq \tau_1$. But then 
$B_{\tau_1} = B_{\tau_0} \in {\partial G_0}$, contrary to the definition of $E$. Therefore, 
$B_{\tau_0}\in H$ and $F$ holds. \qed
\medskip

We also shall need the following result (see \cite{Conway}, Proposition~21.1.17).

\begin{proposition} [Continuity principle] \label{atomless}
If $G$ is a hyperbolic open set and $a \in G$, then the harmonic measure $\omega_{G}(a, \,.\,)$ is atomless.
\end{proposition}

\noindent{\bf Proof of Theorem~\ref{surjective}.} It will be enough to construct a $2$-valent mapping 
$\phi \colon \D \to \D$ such that $\varphi(\D) = \D^*$ and $C_\varphi \colon H^\Psi \to H^\Psi$ is 
compact. We can then modify $\varphi$ by the same trick as the one used by B. McCluer and J. Shapiro.
Note that every point in $\D^\ast$ is the image by $\e^{- z}$ of two distinct points of $R$, except those 
which are the image of points of the hyperbola $y = (1/x) + 2\pi$, which have only one pre-image.\par 
\smallskip

For a positive integer $n$, set: 
\begin{equation}
b_n = \frac{1}{4n\pi} \, \raise 1pt \hbox{,} 
\end{equation}
and let $\varepsilon_n > 0$ such that:
\begin{equation}\label{inegalite sur Psi}
\frac{\Psi^{- 1} (2/b_{n+1})}{\Psi^{- 1} (1/\varepsilon_n)} \leq \frac{1}{n} \,\cdot
\end{equation}

We now modify the domain $R$, including ``barriers'' in it (not in the sense of potential theory, nor of Perron!) 
in the following way. \par

Let, for every $n \geq 1$, $M_n$ be the intersection point of the horizontal line $y = 4\pi n$ and of the 
hyperbola $y = (1/x) + 2\pi$, that is $M_n = \frac{1}{4 \pi n - 2\pi} + 4 \pi n i$.

Define inductively closed sets $P_{n}^+$ and  $P_{n}^-$, which are like small points of swords (two segments  
and a piece of hyperbola), in the following way: 
\begin{itemize}
\item The lower part of  $P_{n}^+$ and $P_{n}^-$ are horizontal segments of altitude $4n\pi$.
\item Those two  horizontal segments are separated by a small open horizontal segment $H_n$ whose 
middle is $M_n$.
\item The upper part of $P_{n}^+$ is a slant segment whose upper extremity $c_{n}^+$ lies on the 
hyperbola $y = 1/ x$.
\item The upper part of $P_{n}^-$ is a slant segment whose upper extremity $c_{n}^-$ lies on the 
hyperbola $y = (1/ x) + 4\pi$.
\item The curvilinear part of $P_{n}^+$ is supported by the hyperbola $y = 1/ x$.
\item The curvilinear part of $P_{n}^-$ is supported by the hyperbola $y = (1/ x) + 4\pi$.
\item One has $4 (n+1) \pi - \Im c_{n}^\pm > 2 \pi$.
\end{itemize}
\goodbreak

\begin{figure}[ht]
\centering
\includegraphics[width=8cm]{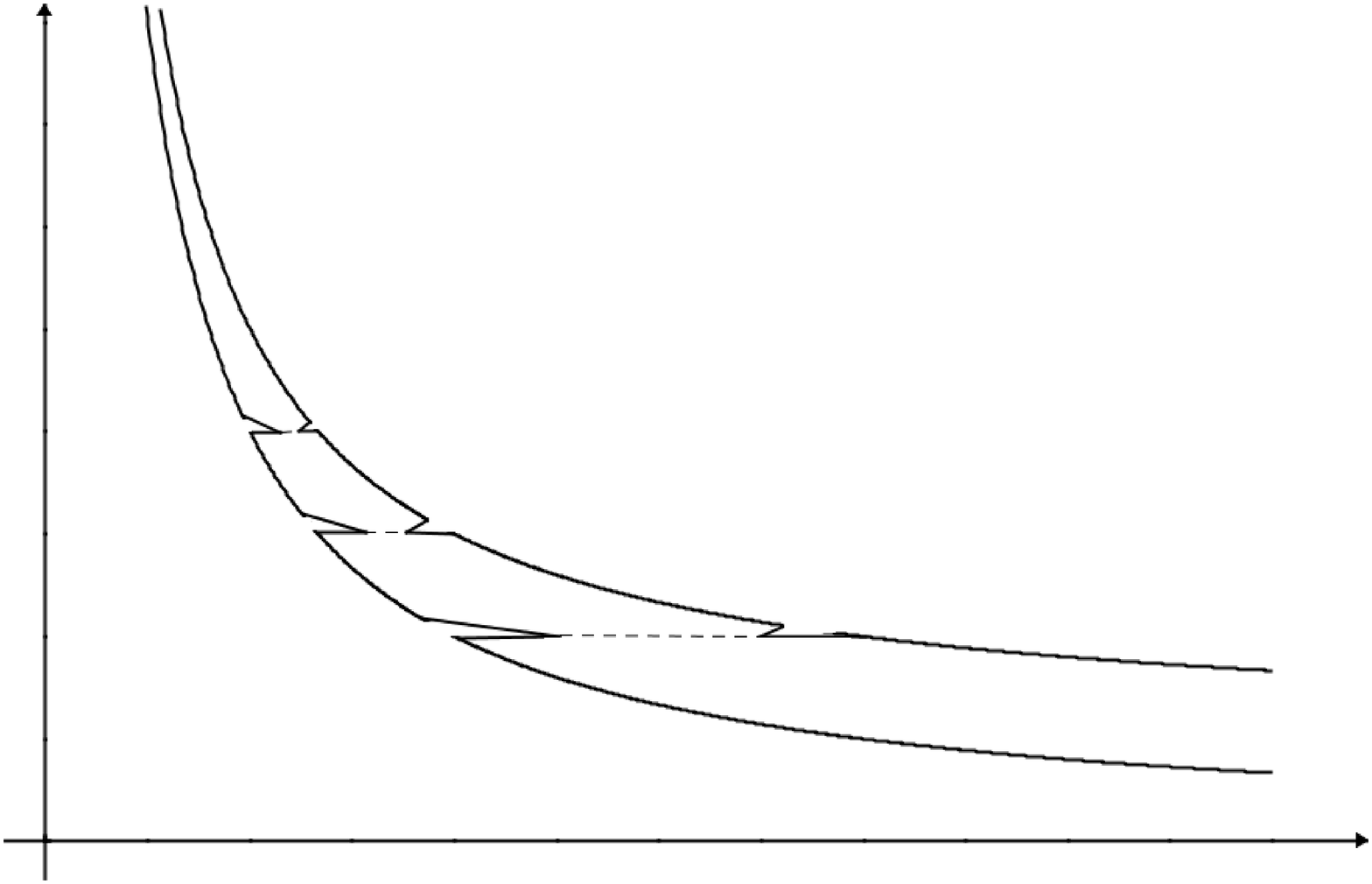}
\end{figure}

The size of the small horizontal holes will be determined inductively in the following way. Fix once and for all  
$a \in R$ such that $\Im a < 4\pi$. Suppose that $H_1, H_2,\ldots, H_{n - 1}$ have already been determined. 
Set:
\begin{equation}
\Omega_n = \big\{z \in R \setminus \bigcup_{j < n} (P_{j}^+ \cup   P_{j}^-) \,;\ \Im z < 4n\pi \big\}.
\end{equation}
We can adjust $H_n$ so small that:
\begin{equation}\label{ajustement H}
\omega_{\Omega_n} (a, H_n)\leq \varepsilon_n. 
\end{equation}
Indeed, $\Omega_n$ is bounded above by the horizontal segment $[b_n + 4in \pi, b_{n - 1} + 4in\pi]$, where 
the point $M_n$ lies. If $H_n = [M_n - \delta, M_n + \delta]$, we see that $H_n$ decreases to the 
singleton $\{M_n\}$ as $\delta$ decreases to zero. Therefore, by Proposition~\ref{atomless}, we can adjust 
$\delta$ so as to realize \eqref{ajustement H}. \par

We now define our modified open set $\Omega$ by the formula 
\begin{equation}
\Omega = R \setminus \bigcup_{n \geq 1} (P_{n}^+ \cup P_{n}^-) = \bigcup _{n \geq 1} \Omega_n.
\end{equation}

It is useful to observe that: 
\begin{equation}\label{borne inf Omega indice n}
\inf_{w \in {\partial \Omega_n}} \Re w = b_n \, .
\end{equation}
This is obvious by the way we defined the upper part of ${\partial}\Omega_n$.\par
\smallskip

Now, we can easily finish the proof. Fix $h\leq b_{1}/ 2$ and let $n$ be the integer such that:
\begin{equation}\label{BN}
b_{n + 1} < 2h \leq b_n \, . 
\end{equation} 
Let $g \colon \D\to \Omega$ be a conformal mapping such that $g(0) = a$. Since $\partial_\infty \Omega$ 
is connected, Caratheodory's Theorem (see \cite{Pomme}) ensures that $g$ can be continuously extended from 
$\overline \D$ onto $\overline \Omega$. More explicitly, using the Moebius transformation 
$T (z) = 1/z$, we see that there exists an automorphism of the extended complex plane such that 
$\overline \Omega$ is sended onto a compact subset of $\C$; so, we can apply to $\Omega$ many results 
stated for bounded domains. For instance, the boundary of $\Omega$ is a continuous path in the extended plane; 
so, by \cite{Conway}, Theorem 14.5.5, $g$ can be extended to a continuous function (for the extended plane 
topology) $g \colon \overline \D \to \overline G$. In particular, $g$ has boundary values $g^\ast$.\par 
We define $\phi = \e^{- g}$.\par
As in the proof of B. McCluer and J. Shapiro (\cite{McCluer-Shapiro}), we have that $\varphi$ is $2$-valent 
(see the remark made at the beginnig of this proof), and we still have $\varphi (\D) = \D^*$, since, in the process 
for constructing $\Omega$ from $R$, for every point of $\D^\ast$, at least one of the preimages by $\e^{- z}$ 
in $R$ has not been removed. Observe that, in particular, we did not remove any point in the hyperbola 
$y = (1/x) + 2\pi$, thanks to the choice of $M_n$.\par 
Moreover, $\Omega$ is a Dirichlet domain (because each component of $\partial \Omega$ has more than one 
point: see the comment after Definition~19.7.1 in \cite{Conway}), so we can use the conformal invariance. Then 
by \eqref{equiv 3}, \eqref{ajustement H}, \eqref{borne inf Omega indice n} and by the hole principle, we see that, 
if $A =\{\Re g^*(\e^{it}) < 2h \}$:  
\begin{align}\label{Carleson}
\rho_{\varphi} (h) 
& \leq m_{\varphi} (\{\vert z \vert > 1 - h\}) =m (\{\e^{-\Re g^*(\e^{it})} > 1 - h\})  \\
& =m (\{\Re g^*(\e^{it}) < \log (1/ 1 - h) \}) \notag \\
& \leq m (\{\Re g^*(\e^{it}) < 2h \}) = \omega_\D (0, A) \notag \\
& =\omega_{g (\D)} \big( g (0), g (A) \big) 
=\omega_{\Omega} (a, \{\Re w < 2h \}) \notag \\
& \leq \omega_{\Omega} (a, \{\Re w \leq b_{n} \}) \notag \\
& \leq \omega_{\Omega} (a, {\partial \Omega} \setminus {\partial \Omega_n}) 
\leq \omega_{\Omega_n} (a, H_n) \leq \varepsilon_n. \notag
\end{align} 

It remains to observe that:
\begin{displaymath}
\Delta (h) = \frac{\Psi^{- 1} (1/ h)}{\Psi^{- 1} (1/ \rho_{\varphi} (h))} 
\leq  \frac{\Psi^{- 1} (2/ b_{n+1})}{\Psi^{- 1} (1/ \varepsilon_n)} 
\leq \frac{1}{n} \leq Ch\,, 
\end{displaymath}
in view of \eqref{inegalite sur Psi} and of the choice of $n$, $C$ being a numerical constant. We should point 
out the fact that we applied the hole principle to the domains $G_0 =\Omega_n$ and $G_1 =\Omega$ 
and that this was licit because the assumptions of the hole principle (in particular the inclusion 
${\partial \Omega_n} \subseteq {\partial \Omega}\cup H_n$) are satisfied. We have therefore proved that: 
\begin{displaymath}
\lim_{h \mathop{\to}\limits^{>} 0} \Delta (h) = 0\,,
\end{displaymath}
and this ends, as we already explained, the first part of the proof of Theorem~\ref{surjective}. \par
\medskip

To prove the last part, let us remark that in \eqref{inegalite sur Psi} we may take $\eps_n$ arbitrarily small. If one 
takes $\eps_n \leq \e^{- n}$, one has, for some constant $c > 0$, $\rho_\phi (h) \leq \e^{- c/h}$, by using 
\eqref{BN} and \eqref{Carleson}. In particular, $\rho_\phi (h) \leq C\, h^\alpha$ for every $\alpha > 1$. By 
Luecking's criterion, that implies that $C_\phi \in S_p (H^2)$ for every $p > 0$ (see \cite{JFA}, Corollary~3.2).
\qed
\goodbreak
\bigskip

\noindent{\bf Remark.} Let us note that our result is stronger than McCluer-Shapiro's, since our $C_\phi$ is in all 
the Schatten classes $S_p (H^2)$, $p >0$. Though our construction follows McCluer-Shapiro's, it is the 
introduction of the ``barriers'' $P_n^+$ and $P_n^-$ which allows to get this improvement.\par\goodbreak


\section{Composition operators with closed range}\label{section closed range}

In \cite{Cima}, J. Cima, J. Thomson and W. Wogen gave a characterization of composition operators 
$C_\phi \colon H^p \to H^p$ with closed range. This characterization involves the Radon-Nikodym 
derivative of the restriction to $\partial \D$ of $m_\phi$. They found it not satisfactory, and asked a  
characterization with the range of $\phi$ itself. N. Zorboska (\cite{Nina}) gave such a characterization, but 
her statement is somewhat complicated. We shall give here more explicit characterizations, either in terms of 
the Nevanlinna counting function $N_\phi$, or in terms of the Carleson measure $m_\phi$. 

\begin{theoreme}\label{closed range Nevanlinna}
Let $\phi \colon \D \to \D$ be a non-constant analytic self map. Then the composition operator 
$C_\phi \colon H^p \to H^p$, $1 \leq p < \infty$, has a closed range if and only if there is a constant 
$c > 0$ such that, for $0 < h < 1$,
\begin{equation}
\qquad \quad \frac{1}{A \big( S (\xi, h) \big)} \int_{S (\xi, h)} N_\phi (z)\,dA (z) \geq c\,h\,, 
\qquad   \forall \xi \in \partial \D\,.
\end{equation}
\end{theoreme}

Theorem~\ref{closed range Nevanlinna} will follow immediately from the next theorem, applied to 
$\mu = m_\phi$, and from \cite{LLQR-N}, Theorem~4.2.

\begin{theoreme}\label{closed range}
Let $\mu$ be a finite positive measure on $\overline \D$. Assume that the canonical map 
$J \colon H^p \to L^p (\mu)$ is continuous, $1 \leq p < \infty$. Then $J$ is one-to-one and has a 
closed range if and only if there is a constant $c > 0$ such that, for $0 < h < 1$,
\begin{equation}\label{mino}
\qquad \qquad \quad \mu \big[ W (\xi, h) \big] \geq c\,h\,, \qquad \qquad  \forall \xi \in \partial \D\,.
\end{equation}
\end{theoreme}

\noindent{\bf Proof.} 1) Assume that $J$ has a closed range. By making a rotation on the variable $z$, 
we only have to find a constant $c > 0$ such that 
\begin{equation}\label{mino S}
\mu (S_h) \geq c\,h \,,
\end{equation}
for $h > 0$ small enough, where $S_h = S (1, h)$.\par
Since $J$ is one-to-one, there is a constant $C > 0$ such that:
\begin{equation}\label{mino f}
\qquad \| f \|_{L^p(\mu)}^p \geq C^p\, \| f \|_p^p \,, \qquad \forall f \in H^p.
\end{equation}

We are going to test \eqref{mino f} on
\begin{equation}
f_N (z) = \bigg( \frac{1 + z}{2} \bigg)^N\,.
\end{equation}

It is classical that there is a constant $c_p > 0$ such that:
\begin{equation}\label{Wallis}
\| f_N \|_p^p = \int_{-\pi}^\pi \Big| \cos \frac{t}{2} \Big|^{pN} \,dt 
\geq \frac{c_p}{\sqrt N} \,\cdot
\end{equation}

Now, since $|z + 1|^2 + |z - 1|^2 = 2 (|z|^2 + 1) \leq 4$ for every $z \in \overline \D$, one has:
\begin{displaymath}
| f_N (z) | \leq \Big( 1 - \frac{|z - 1|^2}{4} \Big)^{N/2} \leq \e^{-\frac{N}{8}\,|z - 1|^2} \,.
\end{displaymath}
Hence, using $| f_N (z) | \leq 1$ when $|z - 1| \leq h$, one has:
\begin{align*}
\| f_N \|_{L^p(\mu)}^p 
& \leq \mu (S_h) + \int_{|z - 1| > h} \e^{- p \frac{N}{8}\,|z - 1|^2} \,d\mu \\
& = \mu (S_h) + \int_0^{\e^{- pN h^2/8}} 
\mu \big( \{\e^{- p \frac{N}{8}\,|z - 1|^2} > u\} \big) \,du \,,
\end{align*}
that is, making the change of variable $u = \e^{- p \frac{N}{8}\,x^2}$,
\begin{displaymath}
\| f_N \|_{L^p(\mu)}^p  \leq 
\mu (S_h) + \int_h^\infty \mu (\{ | z - 1| \leq x\} )\,
\frac{pN}{4}\,x\, \e^{- p \frac{N}{8}\,x^2}\, dx \,.
\end{displaymath}

Now, the continuity of $J$ means, by Carleson's Theorem see \cite{Du}, Theorem~9.3), that there is a constant 
$K > 0$ such that:
\begin{equation}
\qquad \quad \sup_{|\xi| = 1} \mu \big( S (\xi, x) \big) \leq K\,x \,, \qquad 0 \leq x < 1\,.
\end{equation}
We get hence:
\begin{align*}
\| f_N \|_{L^p(\mu)}^p  
& \leq \mu (S_h) + \int_h^\infty K\,x\, \frac{pN}{4}\,x\, \e^{- p \frac{N}{8}\,x^2}\, dx \\
& = \mu (S_h) + \frac{K \sqrt{8}}{\sqrt p}\, \frac{1}{\sqrt N} 
\int_{h \sqrt{\frac{pN}{8}}}^\infty  y^2 \,\e^{-y^2} \,dy \,.
\end{align*}
We take now for $N$ the smaller integer $> 1/ h^2$, multiplied by some constant integer $a_p$, large enough 
to have:
\begin{displaymath}
\frac{K \sqrt{8}}{\sqrt p} 
\int_{\sqrt{\frac{p\,a_p}{8}}}^\infty  \quad y^2 \,\e^{-y^2} \,dy  \leq \frac{c_p\, C^p}{2} \,\cdot
\end{displaymath}
We get then, from \eqref{mino f} and \eqref{Wallis}:
\begin{displaymath}
\mu (S_h) \geq \frac{C^p \,c_p}{2}\, \frac{1}{\sqrt N} \, \raise 1pt \hbox{,}
\end{displaymath}
which gives \eqref{mino S}.
\par\medskip

2) Conversely, assume that \eqref{mino} holds. Since the disk algebra $A (\D)$ is dense in $H^p$, it suffices to 
show that there exists a constant $C > 0$ such that $\| f \|_{L^p (\mu)} \geq C\, \| f \|_p$ for every 
$f \in A (\D)$.\par
Let $f \in A (\D)$ such that $\| f \|_p = 1$. Choose an integer $N$ such that:
\begin{displaymath}
\frac{1}{N} \sum_{n = 1}^N | f (\e^{2\pi i n/N}) |^p 
\geq \frac{1}{2} \int_{\partial \D} | f (\xi) |^p\,dm (\xi) = \frac{1}{2}\,\raise 1pt \hbox{,}
\end{displaymath}
and such that, due to the uniform continuity of $f$,
\begin{displaymath}
z, z' \in \overline \D \quad \text{and} \quad | z - z' |\leq \frac{2\pi}{N} \qquad 
\Longrightarrow \qquad |f (z) - f (z') | \leq \frac{1}{2^{(p + 1)/p}}\, \cdot
\end{displaymath}
Then, setting $W_n = W (\e^{2\pi i n/N}, \pi /N)$, $1 \leq n \leq N$, one has:
\begin{displaymath}
\| f \|_{L^p(\mu)}^p  = \int_{\overline \D} | f |^p\,d\mu 
\geq \sum_{n = 1}^N \int_{W_n} | f |^p\,d\mu\,.
\end{displaymath}
If we choose $z_n \in W_n$ such that $| f (z_n)| = \min_{z \in W_n} | f(z)|$, we get, using \eqref{mino}:
\begin{displaymath}
\| f \|_{L^p(\mu)}^p 
\geq \sum_{n = 1}^N  | f (z_n)|^p \, \mu (W_n)  
\geq \frac{c \pi}{N} \sum_{n = 1}^N  | f (z_n)|^p \,.
\end{displaymath}
Since $A^p \leq 2^{p - 1} [ (A - B)^p + B^p]$, by H\"older's inequality, one has:
\begin{displaymath}
| f (z_n)|^p \geq \frac{1}{2^{p - 1}}\, |f ( \e^{2 \pi i n /N}) |^p - | f (z_n) - f (\e^{2\pi in /N} ) |^p  
\end{displaymath}
and hence:
\begin{displaymath}
\| f \|_{L^p(\mu)}^p 
\geq \frac{c \pi}{N} \sum_{n = 1}^N 
\bigg[\frac{1}{2^{p - 1}}\, |f ( \e^{2 \pi i n /N}) |^p - | f (z_n) - f (\e^{2\pi in /N} ) |^p \bigg]\,.
\end{displaymath}
Now, since $z_n \in W_n$, one has:
\begin{displaymath}
|z_n - \e^{2\pi in /N} | 
\leq \bigg| z_n - \frac{z_n}{|z_n|} \bigg| + \bigg|\frac{z_n}{|z_n|} - \e^{2\pi in /N} \bigg| 
\leq \frac{\pi}{N} + \frac{\pi}{N} = \frac{2\pi}{N}\,; 
\end{displaymath}
therefore $| f (z_n) - f (\e^{2\pi in /N} ) | \leq 1/2^{p + 1}$ and we get:
\begin{align*}
\| f \|_{L^p(\mu)}^p 
& \geq c \pi \,\bigg[\frac{1}{N} \sum_{n = 1}^N  
\frac{1}{2^{p - 1}} \, | f ( \e^{2 \pi i n /N}) |^p - \frac{1}{2^{p + 1}} \bigg] \\
& \geq c \pi \Big( \frac{1}{2^{p - 1}}\,\frac{1}{2}  - \frac{1}{2^{p + 1}} \Big) 
= \frac{c \pi}{2^{p + 1}}\,\cdot
\end{align*}
That ends the proof of Theorem~\ref{closed range}. \qed
\medskip

\noindent{\bf Remark.} To make the link with Cima-Thomson-Wogen's criterion, we shall see that 
condition~\ref{mino} implies that the restriction of $\mu$ to the boundary $\T = \partial \D$ of the disk 
dominates the Lebesgue measure $m$. In fact, let $I$ be an arc of $\T$. If $m (I) = h$, we can write:
\begin{displaymath}
I = \bigcap_{n \geq 1} \bigcup_{j = 1}^n W (\xi_{n, j}, h/2n) \,,
\end{displaymath}
with disjoint windows $W (\xi_{n, 1}, h/2n), \ldots, W (\xi_{n, n}, h/2n)$; hence:
\begin{displaymath}
\mu (I) = \lim_{n \to \infty} \sum_{j=1}^n \mu [W (\xi_{n, j}, h/2n)] \geq c \sum_{j=1}^n \frac{h}{2n} 
= \frac{c}{2}\, h\,.
\end{displaymath}


\section{Composition operators in Schatten classes}\label{Schatten}

In \cite{Luecking}, D. Luecking characterized composition operators $C_\phi \colon H^2 \to H^2$ which 
are in the Schatten classes, by using, essentially, the $m_\phi$-measure of Carleson windows. Five years 
later, D. Luecking and K. Zhu (\cite{Luecking-Zhu}) characterized them by using the Nevanlinna counting 
function of $\phi$. We shall see in this section how the result of \cite{LLQR-N} makes these two 
characterizations directly equivalent. 
\par
It will be convenient here to work with \emph{modified} Carleson windows, namely:
\begin{displaymath}
W_{n, j} = \bigg\{ z\in \overline{\D}\,; \ 
1 - 2^{- n} \leq | z | \leq 1 \text{ and }  
\frac{(2j - 1)\pi}{2^n} \leq \arg z < \frac{(2j + 1 )\pi}{2^n} \bigg\} 
\end{displaymath}
($j = 0, 1, \ldots, 2^n - 1$, $n = 1, 2, \ldots$). We shall say that $W_{n, j}$ is the Carleson window  
centered at $\e^{2\pi i j/2^n}$ with size $2^{- n}$.

\begin{theoreme}\label{equiv Luecking}
For $p >0$ the two following conditions are equivalent:\par
\smallskip
$a)$ $\displaystyle \frac{N_\phi (z)}{\log (1/ |z|)} \in L^{p/2} (\lambda)$\,, where 
$d\lambda (z) = (1 - |z|)^{-2}\,dA (z)$ and $A$ is the normalized area measure on $\D$;\par
\smallskip
$b)$  
$\displaystyle \sum_{n = 1}^\infty \sum_{j = 0}^{2^n - 1} \big[ 2^n \, m_\phi (W_{n, j}) \big]^{p/2} < \infty$\,.
\end{theoreme} 

Condition $b)$ in the last theorem yields that $\lim_{n\to\infty}\max_j 2^n m_\varphi (W_{n,j}) = 0$, and 
it is not difficult to see that this implies that $m_\varphi (\partial \D) = 0$, or equivalently, that 
$|\varphi^\ast|<1$ almost evereywhere on $\partial \D$. In this situation we know 
(\cite{JFA}, Proposition~3.3) that $b)$ in Theorem~\ref{equiv Luecking} is equivalent to Luecking's condition 
in \cite{Luecking}. In fact the characterization of belonging to a Schatten class in \cite{Luecking} includes the 
requirement $m_\varphi(\partial \D) = 0$.\par
\medskip

\noindent{\bf Proof.} We may, and do, assume that $\phi (0) = 0$.\par
\smallskip

1) Assume first that condition~$b)$ is satisfied. Let:
\begin{displaymath}
R_{n, j} = \Big\{z\in \D\,;\ 1 - 2^{-n} \leq |z| < 1 - 2^{-n - 1}\ \text{and}\  
\frac{(2j - 1) \pi}{2^n} \leq \arg z < \frac{(2j+ 1) \pi}{2^n}\,\Big\} 
\end{displaymath}
be the (disjoint) Luecking windows ($0 \leq j \leq 2^n - 1$, $n \geq 0$). One has 
$R_{n, j} \subseteq W_{n, j}$.\par
By \cite{LLQR-N}, Theorem~3.1, there are a constant $C > 0$ and an integer $K$ 
such that $N_\phi (z) \leq C\,m_\phi (\widetilde W_{n, j})$, for every $z \in R_{n, j}$, where 
$\widetilde W_{n, j}$ is the window centered at $\e^{2\pi i j/2^n}$, as $W_{n, j}$, but with size 
$2^{K - n}$. The windows $W_{n - K, j}$, $j = 0, 1, \ldots, 2^{n - K} - 1$, have the same size as the windows 
$\widetilde W_{n, j}$, but may have a different center; nevertheless, each $\widetilde W_{n, j}$ can be 
covered with two windows $W_{n - K, l}$: for $n > K$,  
$\widetilde W_{n, j} \subseteq W_{n - K, l} \cup W_{n - K, l + 1}$, for some $l = 1, 2, \ldots, 2^{n - K}$ 
(where $l + 1$ is understood as $0$ if $l = 2^{n - K} - 1$), we get (we shall use $\lesssim$ to mean $\leq$ 
up to a constant):
\begin{align*}
\int_\D \frac{\big(N_\phi (z) \big)^{p/2}}{(1 - |z|)^{\frac{p}{2} + 2}}\,dA (z) 
& \leq \sum_{n, j} \int_{R_{n, j}} (2^n)^{\frac{p}{2} + 2} \big( N_\phi (z) \big)^{p/2}\,dA (z) \\
& \lesssim \sum_{n, j} \int_{R_{n, j}} 
(2^n)^{\frac{p}{2} + 2} \big( m_\phi (\widetilde W_{n, j}) \big)^{p/2}\,dA (z) \\
& \lesssim \sum_{n, j} (2^n)^{p/2}\, \big(m_\phi (\widetilde W_{n, j}) \big)^{p/2} \\
& \lesssim \sum_{\nu, l} (2^\nu)^{p/2}\, \big(m_\phi (W_{\nu, l}) \big)^{p/2} < \infty,
\end{align*}
and $a)$ holds.\par\medskip

2) Conversely, assume that $a)$ is satisfied. We shall use the following inequality, whose proof will be 
postponed (for $p \geq 2$, \eqref{inegalite utile} follows directly from \cite{LLQR-N}, Theorem~4.2,  and 
H\"older's inequality):
\begin{equation}\label{inegalite utile}
[m_\phi (W_{n, j})]^{p/2} 
\lesssim \frac{1}{A (\widetilde W_{n, j})} \int_{\widetilde W_{n, j}} [N_\phi (z)]^{p/2}\,dA (z)\,,
\end{equation}
where $\widetilde W_{n, j}$ is a window with the same center as $W_{n, j}$ but with a bigger proportional 
size; say of size $2^{- n + L}$. We get:
\begin{align*}
\sum_{n, j} [2^n \, m_\phi (W_{n, j} )]^{p/2} 
& \lesssim \sum_{n, j} 2^{np/2}\, 2^{2n}\, \int_{\widetilde W_{n, j}} [N_\phi (z)]^{p/2}\, dA (z) \\
& = \int_{\D} \bigg( \sum_n 2^{n (2 + \frac{p}{2})} \Big[\sum_j \ind_{\widetilde W_{n, j}} (z) \Big] \bigg)
\, [N_\phi (z)]^{p/2}\, dA (z) \,.
\end{align*}
Let $k = 0, 1, \ldots$ such that $1 - 2^{- k + 1} < |z| \leq 1 - 2^{- k}$. One has $z \in \widetilde W_{n, j}$ only 
if $n \leq k + L$, and then, for each such $n$, $z$ is at most in $2^L$ windows $\widetilde W_{n, j}$. It 
follows that:
\begin{displaymath}
\sum_n 2^{n (2 + \frac{p}{2})} \sum_j \ind_{\widetilde W_{n, j}} (z) 
\leq 2^{(k + L + 1)(2 + \frac{p}{2})} \times 2^L \,.
\end{displaymath}
But $|z| \geq 1 - 2^{- k + 1}$ implies $2^{(k + L  + 1)(2 + \frac{p}{2})} \leq C_p/ (1 - |z|)^{2 + \frac{p}{2}}$;  
hence:
\begin{align*}
\sum_{n, j} [2^n \, m_\phi (W_{n, j} )]^{p/2} 
\lesssim \int_{\D} \frac{[N_\phi (z)]^{p/2}}{(1 - |z|)^{\frac{p}{2} + 2}}\,dA (z) < \infty\,,
\end{align*}
and $b)$ holds. \par
\medskip\goodbreak

It remains to show \eqref{inegalite utile}.\par
By \cite{LLQR-N}, Theorem~4.1, we can find a window $W$ with the same center as $W_{n, j}$, but with greater 
size $ch$ ($h = 2^{-n}$ is the size of the window $W_{n, j}$), such that:
\begin{displaymath}
m_\phi ( W_{n, j}) \lesssim \sup_{w \in W} N_\phi (w).
\end{displaymath}
There is hence some $w_0 \in W$ such that:
\begin{displaymath}
m_\phi ( W_{n, j}) \lesssim N_\phi (w_0).
\end{displaymath}

Take $R = |w_0| + ch$ (one has $R \geq 1$ since $w_0 \in W$ and $W$ has size $c h$) and set 
$\phi_0 (z) = \phi (z) /R$. One has $N_{\phi_0} (z) = N_\phi (Rz)$ for $|z| < 1/R$ and $N_{\phi_0} (z) = 0$ if 
$|z| \geq 1/R$.\par
Let now $u$ be the upper subharmonic regularization of $N_{\phi_0}$ (\cite{Luecking-Zhu}, Lemma~1, 
and its proof page 1140): $u$ is a subharmonic function on $\D\setminus\{0\}$ such that $u \geq N_{\phi_0}$ 
and $u =  N_{\phi_0}$ almost everywhere, with respect to $dA$.
\par\smallskip

A result of C. Fefferman and E. M. Stein (\cite{Fefferman-Stein}, Lemma~2), generously attributed by them 
to Hardy and Littlewood, asserts that for any $q > 0$, there exists a constant $C = C (q)$ such that
\begin{equation}\label{Feff-Stein}
[u (a)]^q \leq \frac{C}{A \big( D (a, r) \big)} \int_{D (a, r)} [u (z)]^q\, dA (z)
\end{equation}
for every nonnegative subharmonic function $u$ on a domain $G$ and every disk $D (a, r) \subseteq G$ 
(see also \cite{Luecking-Zhu}, Lemma~3).\par\smallskip

If $\Delta$ is the disk centered at $w_0/R$ and of radius $1 - |w_0|/R$ (which is contained in 
$\D \setminus\{0\}$ since $R > |w_0|$), one has, by \eqref{Feff-Stein}:
\begin{align*}
[N_{\phi} (w_0)]^{p/2} 
& = [N_{\phi_0} (w_0/ R)]^{p/2} \leq [ u (w_0/ R)]^{p/2} \\
& \leq \frac{C}{A (\Delta)} \int_\Delta [u (z)]^{p/2} \,dA (z) \\ 
& = \frac{C}{A (\Delta)} \int_\Delta [N_{\phi_0} (z)]^{p/2} \,dA (z) \\
& = \frac{C}{A (\Delta)} \int_{\Delta \cap D (0, 1/R)} [N_{\phi} (Rz)]^{p/2} \,dA (z) \\
& = \frac{C}{A (\tilde \Delta)} \int_{\tilde \Delta \cap \D} [N_\phi (w)]^{p/2} \,dA (w)\,,
\end{align*}
where $\tilde \Delta = D (w_0, R - |w_0|) = D (w_0, ch)$. \par
Since the center $w_0$ of  $\tilde \Delta$ is in $\D$, $\tilde \Delta \cap \D$ contains more than a quarter of 
$\tilde \Delta$ (at least for $c h \leq 1$), and hence 
$A (\tilde \Delta \cap \D) \geq A (\tilde \Delta) /4 = c^2 h^2/4\pi$. Now, let $\tilde W_{n, j}$ be the 
window with the same center as $W_{n, j}$ and of size $2c h$. Since $2c h \geq c h + (1 - |w_0|)$,  
$\tilde W_{n, j}$ contains $\tilde \Delta \cap \D$ and 
$A (\tilde W_{n, j}) \approx h^2 \approx A (\tilde \Delta)$ 
($\approx$ meaning that the ratio is between two absolute constants). We therefore get: 
\begin{displaymath}
[N_{\phi} (w_0)]^{p/2} 
\lesssim \frac{1}{A (\tilde W_{n, j})} \int_{\tilde W_{n, j}} [N_\phi (w)]^{p/2} \,dA (w)\,,
\end{displaymath}
proving \eqref{inegalite utile}.\qed


\bigskip

\vbox{\noindent{\it 
{\rm Pascal Lef\`evre}, Univ Lille Nord de France F-59\kern 1mm 000 LILLE, FRANCE\\
UArtois, Laboratoire de Math\'ematiques de Lens EA~2462, \\
F\'ed\'eration CNRS Nord-Pas-de-Calais FR~2956, \\
F-62\kern 1mm 300 LENS, FRANCE \\ 
pascal.lefevre@euler.univ-artois.fr 
\smallskip

\noindent
{\rm Daniel Li}, Univ Lille Nord de France F-59\kern 1mm 000 LILLE, FRANCE\\
UArtois, Laboratoire de Math\'ematiques de Lens EA~2462, \\
F\'ed\'eration CNRS Nord-Pas-de-Calais FR~2956, \\
Facult\'e des Sciences Jean Perrin,\\
Rue Jean Souvraz, S.P.\kern 1mm 18, \\
F-62\kern 1mm 300 LENS, FRANCE \\ 
daniel.li@euler.univ-artois.fr
\smallskip

\noindent
{\rm Herv\'e Queff\'elec}, Univ Lille Nord de France F-59\kern 1mm 000 LILLE, FRANCE\\
USTL, Laboratoire Paul Painlev\'e U.M.R. CNRS 8524, \\
F-59\kern 1mm 655 VILLENEUVE D'ASCQ Cedex, FRANCE \\ 
queff@math.univ-lille1.fr
\smallskip

\noindent
{\rm Luis Rodr{\'\i}guez-Piazza}, Universidad de Sevilla, \\
Facultad de Matem\'aticas, Departamento de An\'alisis Matem\'atico,\\ 
Apartado de Correos 1160,\\
41\kern 1mm 080 SEVILLA, SPAIN \\ 
piazza@us.es\par}
}

\end{document}